\documentclass[oneside, reqno]{amsart}

\usepackage{xypic}
\input xy
\xyoption{all}
\usepackage{epsfig}
\usepackage{amsthm}
\usepackage{amssymb}
\usepackage{amsmath}
\usepackage{amscd}
\usepackage{color}
\usepackage[T1]{fontenc}
\usepackage{stackengine}
\usepackage{upgreek}
\usepackage[font=scriptsize]{caption}
\usepackage{wrapfig}
\stackMath

\newcommand{\bg}{\begin{equation}}
\newcommand{\ed}{\end{equation}}
\newcommand{\bga}{\begin{eqnarray}}
\newcommand{\eda}{\end{eqnarray}}
\newcommand{\pf}{\textbf{Proof:\ }}

\def\cbdu{\par{\raggedleft$\Box$\par}}

\newtheorem {Theorem}  {Theorem}

\numberwithin{Theorem}{section}

\newtheorem {Lemma}[Theorem]  {Lemma}

\theoremstyle{definition}
\newtheorem{Definition}[Theorem]{Definition}
\theoremstyle{remark}
\newtheorem{Remark}[Theorem]{\bf Remark}

\expandafter\chardef\csname pre amssym.def
at\endcsname=\the\catcode`\@ \catcode`\@=11
\def\undefine#1{\let#1\undefined}
\def\newsymbol#1#2#3#4#5{\let\next@\relax
 \ifnum#2=\@ne\let\next@\msafam@\else
 \ifnum#2=\tw@\let\next@\msbfam@\fi\fi
 \mathchardef#1="#3\next@#4#5}
\def\mathhexbox@#1#2#3{\relax
 \ifmmode\mathpalette{}{\m@th\mathchar"#1#2#3}%
 \else\leavevmode\hbox{$\m@th\mathchar"#1#2#3$}\fi}
\def\hexnumber@#1{\ifcase#1 0\or 1\or 2\or 3\or 4\or 5\or 6\or 7\or 8\or
 9\or A\or B\or C\or D\or E\or F\fi}

\font\teneufm=eufm10 \font\seveneufm=eufm7 \font\fiveeufm=eufm5
\newfam\eufmfam
\textfont\eufmfam=\teneufm \scriptfont\eufmfam=\seveneufm
\scriptscriptfont\eufmfam=\fiveeufm

\catcode`\@=\csname pre amssym.def at\endcsname

\newcounter{remark}
\def  \12  {{\frac{1}{2}}}

\def\build#1_#2^#3{\mathrel{\mathop{\kern 0pt#1}\limits_{#2}^{#3}}}

\numberwithin{equation}{section}

\begin{document}

\title[Uniqueness and non-uniqueness for dyadic MHD]{Uniqueness and non-uniqueness results for dyadic MHD models}


\author [Mimi Dai]{Mimi Dai}

\address{Department of Mathematics, Statistics and Computer Science, University of Illinois at Chicago, Chicago, IL 60607, USA}
\email{mdai@uic.edu} 

\author [Susan Friedlander]{Susan Friedlander}

\address{Department of Mathematics, University of Southern California, Los Angeles, CA 90089, USA}
\email{susanfri@usc.edu}


\begin{abstract}

We construct non-unique Leray-Hopf solutions for some dyadic models for magnetohydrodynamics when the intermittency dimension $\delta$ is less than 1. In contrast, uniqueness of Leray-Hopf solution is established in the case of $\delta\geq 1$. 
Analogous results on uniqueness and non-uniqueness of Leray-Hopf solution are also obtained for dyadic models of MHD with fractional diffusion. 

\bigskip


\hspace{0.02cm}CLASSIFICATION CODE: 35Q35, 76D03, 76W05.
\end{abstract}

\maketitle

\section{Introduction}


\subsection{Magnetohydrodynamics}
In geophysics and astrophysics, incompressible magnetohydrodynamics (MHD) governed by the equations
\begin{subequations}
\begin{align}
u_t+(u\cdot\nabla) u-(B\cdot\nabla) B+\nabla P=&\ \nu\Delta u+f, \label{mhda}\\
B_t+(u\cdot\nabla) B-(B\cdot\nabla) u =&\ \mu\Delta B, \label{mhdb}\\
\nabla \cdot u= 0, \ \ \nabla \cdot B=&\ 0, \label{mhdc}
\end{align}
\end{subequations}
is a fundamental model in the investigation of electrically conducting fluids. 
In this system
the vector fields $u$ and $B$ represent the fluid velocity and magnetic field, respectively; 
the scalar function $P$ denotes the pressure; the parameters $\nu$ and $\mu$ denote respectively the viscosity and the magnetic resistivity; and $f$ stands for an external force acting on the fluid.
When $B=0$, system (\ref{mhda})-(\ref{mhdc}) becomes the Navier-Stokes equation (NSE) (\ref{nse}) which will be discussed later. 

It is evident that the MHD system inherits challenges from the NSE, but also exhibits its own complexity which is mainly caused by the nonlinear interactions between the fluid velocity field and the magnetic field. The unsolved problems for the NSE usually also hang in the air for the MHD system. In particular, it is not clear whether either the NSE or MHD has a classical solution for all the time, given arbitrary initial data.  The concept of Leray-Hopf solution for the NSE was introduced by Leray \cite{Le} and Hopf \cite{Ho}. A Leray-Hopf solution is a weak solution in the standard distributional sense, which satisfies the basic energy inequality. Such a concept is naturally adapted to other partial differential equations. Since the pioneering work of Leray, the well-posedness problem for the Leray-Hopf solutions to the NSE in three dimensional (3D) space is still open. In particular, the uniqueness of a Leray-Hopf solution to the 3D NSE remains unsolved. This is the case for the 3D MHD (\ref{mhda})-(\ref{mhdc}) as well.
Nonetheless, wild weak solutions have been constructed for the ideal MHD, i.e. $\nu=\mu=0$ and $f=0$ in (\ref{mhda})-(\ref{mhdc}) by Beekie, Buckmaster and Vicol \cite{BBV}. The weak solutions constructed in \cite{BBV} have finite total energy, but do not conserve the magnetic helicity which is an invariant quantity for smooth solutions.

The main objective of this paper is to investigate the problem of uniqueness of Leray-Hopf solutions for the diffusive dyadic models of the MHD system (\ref{mhda})-(\ref{mhdc}). 
The following dyadic model for the MHD system was proposed in \cite{Dai-20},
\begin{subequations}
\begin{align}
\frac{d}{dt}a_j+\nu\lambda_j^2 a_j
&+\kappa_1\left(\lambda_j^{\theta}a_ja_{j+1}-\lambda_{j-1}^{\theta}a_{j-1}^2\right) \notag \\
&-\kappa_2\left(\lambda_{j}^{\theta}b_jb_{j+1}-\lambda_{j-1}^{\theta}b_{j-1}^2\right)=f_j, \label{gen1}\\
\frac{d}{dt}b_j+\mu\lambda_j^2 b_j
&+\kappa_2\left(\lambda_j^{\theta}a_jb_{j+1}-\lambda_{j}^{\theta}b_ja_{j+1}\right) =0 \label{gen2}
\end{align}
\end{subequations}
for $j\geq 0$, $\lambda_j=\lambda^j$ with a fixed constant $\lambda>1$, 
and $a_{-1}=b_{-1}=0$. 
The variables $(a_j, b_j)$ are quantities related to the energy in the sense that 
$\frac12a_j^2$ and $\frac12b_j^2$ are the kinetic energy and magnetic energy in the $j$-th shell, respectively. 
The parameter $\theta$ is defined as $\theta=\frac{5-\delta}{2}$, where $\delta\in[0,3]$ is the intermittency dimension for the 3D turbulent vector field (cf. \cite{CD-Kol}). Naturally, $\theta\in[1, \frac52]$. Notice that smaller $\delta$ corresponds to larger $\theta$, and hence stronger nonlinearity.   
The parameters $\kappa_1$ and $\kappa_2$ placed in front of the nonlinear terms
represent the energy transfer direction and strength among shells. 
Similar dyadic models have been presented by physicists for the MHD system, for instance, see \cite{GLPG, PSF}.

Denote the total energy by 
\begin{equation}\notag
E(t)=\frac12\sum_{j=0}^{\infty}\left(a_j^2+b_j^2\right)
\end{equation}
and the flux through the $j$-th shell by 
\begin{equation}\notag
\Pi_j=\lambda_j^\theta (\kappa_1 a_j^2-\kappa_2b_j^2)a_{j+1}, \ \ \ j\geq0.
\end{equation}
The energy balance for the $j$-th shell of the system (\ref{gen1})-(\ref{gen2}) is 
\begin{equation}\notag
\frac{d}{dt}\left(a_j^2+b_j^2\right)=-\nu \lambda_j^2a_j^2-\mu \lambda_j^2b_j^2+\Pi_{j-1}-\Pi_j+f_ja_j.
\end{equation}
Thus system (\ref{gen1})-(\ref{gen2}) obeys the formal energy law
\begin{equation}\notag
\frac{d}{dt}E(t)=-\nu\sum_{j=0}^{\infty} \lambda_j^2a_j^2-\mu\sum_{j=0}^{\infty} \lambda_j^2b_j^2+\sum_{j=0}^{\infty} f_ja_j. 
\end{equation}
It is clear to see that the energy is invariant for (\ref{gen1})-(\ref{gen2}) if $\nu=\mu=0$ and $f_j=0$ for $j\geq 0$. 
We will consider the four particular cases of the general model (\ref{gen1})-(\ref{gen2}) with $\kappa_1=\pm 1$ and $\kappa_2=\pm 1$.


We will provide a definition of Leray-Hopf solutions for dyadic models in analogy with the Leray-Hopf solution for the original MHD equations (\ref{mhda})-(\ref{mhdc}). The main goal is to: (i) establish global in time existence of Leray-Hopf solutions for the dyadic models; (ii) show the uniqueness of Leray-Hopf solution when $\theta\leq 2$; (iii) construct non-unique Leray-Hopf solutions in the case of $\theta>2$. Philosophically, the process of constructing non-uniqueness resembles the convex integration method in the sense that it takes advantage of the forcing term in the construction. Technically it is much simpler than convex integration since no iteration or approximation is involved. 


\medskip

\subsection{Main results for dyadic MHD models}
In this part, we lay out the results regarding Leray-Hopf solutions for the dyadic MHD model (\ref{gen1})-(\ref{gen2}). First, for any initial data with finite total energy, we show the existence of global Leray-Hopf solutions. 


\begin{Theorem}\label{thm-existence}
Let $\theta>0$, $a^0=\{a_j^0\}_{j\geq0}\in l^2$ and $b^0=\{b_j^0\}_{j\geq0}\in l^2$. For any $T>0$, assume 
\[\sum_{j=0}^\infty \lambda_j^{-2}\int_0^T f_j^2(t)\, dt<\infty,\]
i.e. $f\in L^2(0,T; H^{-1})$.
Then there exists a Leray-Hopf solution to system (\ref{gen1})-(\ref{gen2}) accompanied with the initial data $(a^0, b^0)$ on $[0,T]$.
\end{Theorem}

The next result concerns the weak-strong type of uniqueness.

\begin{Theorem}\label{thm-weak-strong}
Let $\theta>0$. Let $(a(t), b(t))$ and $(u(t), v(t))$ be Leray-Hopf solutions to (\ref{gen1})-(\ref{gen2}) with the same initial data $(a^0, b^0)\in l^2\times l^2$. Assume in addition that there is a number $J$ such that 
\begin{equation}\label{assump-a}
|a_j(t)|\leq C_0 \lambda_j^{2-\theta}, \ \ |b_j(t)|\leq C_0 \lambda_j^{2-\theta} \ \ \mbox{for} \ \ j\geq J, \ \ t\in[0,T]
\end{equation}
with a constant $C_0$ depending on $\lambda$ and $\theta$. Then 
\[a_j\equiv u_j, \ \ \ b_j\equiv v_j, \ \ \ \mbox{on} \ \ [0,T] \ \ \mbox{for all} \ \ j\geq 0.\]
\end{Theorem}

As a consequence of Theorem \ref{thm-weak-strong}, the uniqueness of the Leray-Hopf solution in the case of $\theta\leq 2$ follows immediately. 

\begin{Theorem}\label{thm-unique1}
Let $0<\theta\leq 2$. Let $a^0=\{a_j^0\}_{j\geq0}\in l^2$, $b^0=\{b_j^0\}_{j\geq0}\in l^2$ and $f\in L^2(0,T; H^{-1})$. Then the Leray-Hopf solution to (\ref{gen1})-(\ref{gen2}) is unique.
\end{Theorem}

When $\theta>2$, we adapt the construction strategy of \cite{FK} for the dyadic NSE and show that the dyadic MHD models have more than one Leray-Hopf solutions. Specifically, we prove:

\begin{Theorem}\label{thm-nonunique1}
Let $\theta>2$. Let $a^0=0$ and $b^0=0$, i.e $a_j^0=b_j^0=0$ for all $j\geq 0$. There exists $T>0$ and functions $\{f_j(t)\}$ satisfying
\[\sum_{j=0}^\infty \lambda_j^{-2} \int_0^T f_j^2(t)\, dt<\infty,\]
such that system (\ref{gen1})-(\ref{gen2}) with initial data $(a^0, b^0)$ has at least two Leray-Hopf solutions $(a(t), b(t))$, one of which has non-vanishing $a(t)$ and $b(t)$ on $[0,T]$. 
\end{Theorem}




\begin{Remark}
The solutions constructed in Theorem \ref{thm-nonunique1} satisfy the energy identity.  
\end{Remark}

\begin{Remark}
We see that the threshold value of $\theta$ that separates the uniqueness and non-uniqueness results is $\theta=2$. Notice that since $\theta=\frac{5-\delta}{2}$, $\theta=2$ corresponds to the intermittency dimension $\delta=1$. In fact, there is evidence that $\delta=1$ is critical for 3D turbulent flows, see \cite{CSreg}.
\end{Remark}

\medskip

\subsection{Weak solutions for dyadic NSE}
The incompressible Navier-Stokes equation
\begin{equation}\label{nse}
\begin{split}
u_t+(u\cdot\nabla)u+\nabla P=&\ \nu\Delta u+f, \\
\nabla\cdot u=&\ 0,
\end{split}
\end{equation}
is a central topic in the study of fluids. In (\ref{nse}), the vector field $u$,  scalar function $P$, parameter $\nu$ and vector valued function $f$ denote the fluid velocity, pressure, viscosity and external forcing, respectively.
Although there has been much progress in the past century concerning fundamental properties of the NSE, many significant questions remain open. Partly for this reason various so called dyadic models have been proposed. One such model for oceanographic turbulence was presented by Desnyanskiy and Novikov \cite{DN} in 1974 and later with motivation from harmonic analysis by Katz and Pavlovi\'c \cite{KP}. This model takes the form

\begin{equation}\label{nse-dyadic}
\frac{d}{dt}a_j +\nu\lambda_j^2 a_j
 + \lambda_j^{\theta}a_ja_{j+1}
 -\lambda_{j-1}^{\theta}a_{j-1}^2 = f_j,
\end{equation}
for $j\geq 1$ and $a_0=0$. 
A crucial property of this particular model is the persistence of positivity, namely that with nonnegative forcing a solution starting from positive initial data remains positive for all time. This attribute of the system (\ref{nse-dyadic}) was essential for the proof of many interesting results, for example see \cite{BFM, BMR, Ch, CF, CFP1, CFP2}. However, as was recently observed by Filonov and Khodunov \cite{FK} the energy cascade in a turbulent fluid is a random process with no physical reason for the conservation of positivity. Hence there is an intrinsic desirability for techniques that do not depend on positivity. In \cite{FK} the authors introduced a novel approach that does not depend on positivity. They proved for (\ref{nse-dyadic}) existence and uniqueness of Leray-Hopf solution with $\theta\leq 2$ and showed that there exist more than one Leray-Hopf solution with $\theta>2$. Specifically, they proved the following theorem:

\begin{Theorem}\label{thm-nse-FK}\cite{FK}
Let $\theta>2$ and $a^0=0$. There exists $T>0$ and functions $f_j(t)$ satisfying $f=\{f_j\}_{j\geq 0}\in L^2(0,T; H^{-1})$ such that the dyadic NSE model (\ref{nse-dyadic}) with initial data $a^0$ has at least two Leray-Hopf solutions. 
\end{Theorem} 

Returning to the dyadic MHD model (\ref{gen1})-(\ref{gen2}) the delicate interactions between the velocity and the magnetic fields preclude the possibility of making sign choice of the parameters that ensure the persistence of positivity. The techniques that we use to prove the results stated in Theorems \ref{thm-weak-strong} - \ref{thm-nonunique1} are motivated by the approach used for the NSE in \cite{FK} which does not depend on positivity. We observe that in this context the complexity of the nonlinear coupling of the two fields is actually a benefit that gives us additional freedom in constructing a scheme used to prove non-uniqueness.

\medskip

\subsection{Dyadic models with fractional Laplacian scaling}

We note that the dyadic MHD equations (\ref{gen1})-(\ref{gen2}) can also be rescaled to
\begin{equation}\label{sys-3}
\begin{split}
\frac{d}{dt}a_j +\nu\lambda_j^{2\alpha} a_j
 &+ \kappa_1(\lambda_j a_ja_{j+1}-\lambda_{j-1}a_{j-1}^2)\\
&-\kappa_2( \lambda_{j} b_jb_{j+1} -\lambda_{j-1}b_{j-1}^2) = f_j,\\
\frac{d}{dt}b_j +\mu\lambda_j^{2\alpha} b_j
&+\kappa_2( \lambda_ja_jb_{j+1}-\lambda_{j}b_ja_{j+1})= 0,
\end{split}
\end{equation}
for $j\geq 1$, $a_0=b_0=0$ and $\alpha=\frac{1}{\theta}$. 

The analogous dyadic model for the fractional MHD with diffusion terms $(-\Delta)^{\alpha} u$ and $(-\Delta)^{\beta} B$ is
\begin{equation}\label{sys-4}
\begin{split}
\frac{d}{dt}a_j +\nu\lambda_j^{2\alpha} a_j
 &+ \kappa_1(\lambda_j a_ja_{j+1}-\lambda_{j-1}a_{j-1}^2)\\
&-\kappa_2 (\lambda_{j} b_jb_{j+1} -\lambda_{j-1}b_{j-1}^2) = f_j,\\
\frac{d}{dt}b_j +\mu\lambda_j^{2\beta} b_j
&+\kappa_2 (\lambda_ja_jb_{j+1}-\lambda_{j}b_ja_{j+1})= 0,
\end{split}
\end{equation}
with $j\geq 1$, $a_0=b_0=0$, and $\alpha>0, \beta>0$. Obviously (\ref{sys-3}) is a special case of (\ref{sys-4}) with $\alpha=\beta$.
With slight modifications of the proof for Theorem \ref{thm-existence}, we can prove that:
\begin{Theorem}\label{thm-existence-gmhd}
Let $\alpha>0$ and $\beta>0$. Let $a^0=\{a_j^0\}_{j\geq0}\in l^2$ and $b^0=\{b_j^0\}_{j\geq0}\in l^2$.  Assume $f\in L^2(0,T; H^{-\alpha})$ for any $T>0$.
Then there exists a Leray-Hopf solution to system (\ref{sys-4}) accompanied with the initial data $(a^0, b^0)$ on $[0,T]$.
\end{Theorem}

In analogy with Theorem \ref{thm-weak-strong}, we can prove the following weak-strong type of uniqueness for a Leray-Hopf solution to (\ref{sys-4}).

\begin{Theorem}\label{thm-weak-strong-gmhd}
Let $\alpha>0$ and $\beta>0$. Let $(a(t), b(t))$ and $(u(t), v(t))$ be Leray-Hopf solutions to (\ref{sys-4}) with the same initial data $(a^0, b^0)\in l^2\times l^2$. Assume in addition that there is a number $J$ such that 
\begin{equation}\label{assump-gmhd}
|a_j(t)|\leq C_0 \left(\lambda_j^{2\alpha-1}+\lambda_j^{2\beta-1}\right), \ \ |b_j(t)|\leq C_0 \lambda_j^{\alpha+\beta-1}
\end{equation}
for all $j\geq J$ and $t\in [0,T]$,
with a constant $C_0$ depending on $\lambda$ and $\theta$. Then 
\[a_j\equiv u_j, \ \ \ b_j\equiv v_j, \ \ \ \mbox{on} \ \ [0,T] \ \ \mbox{for all} \ \ j\geq 0.\]
\end{Theorem}

The following uniqueness of a Leray-Hopf solution to (\ref{sys-4}) with $\alpha\geq \frac12$ and $\beta\geq \frac12$ is an immediate consequence of Theorem \ref{thm-weak-strong-gmhd}.

\begin{Theorem}\label{thm-unique-gmhd}
Let $\alpha\geq \frac12$ and $\beta\geq \frac12$. Let $a^0=\{a_j^0\}_{j\geq0}\in l^2$, $b^0=\{b_j^0\}_{j\geq0}\in l^2$ and $f\in L^2(0,T; H^{-\alpha})$. Then the Leray-Hopf solution to (\ref{sys-4}) is unique.
\end{Theorem}

We also construct non-unique Leray-Hopf solutions to (\ref{sys-4}) for appropriate values of $\alpha$ and $\beta$. Namely, we will show:

\begin{Theorem}\label{thm-nonunique2}
Let $0<\alpha\leq\beta<\frac12$ and $3\beta-\alpha<1$. Let $a^0=0$ and $b^0=0$. There exists $T>0$ and functions $\{f_j(t)\}$ satisfying
\[\sum_{j=0}^\infty \lambda_j^{-2\alpha} \int_0^T f_j^2(t)\, dt<\infty,\]
such that system (\ref{sys-4}) with initial data $(a^0, b^0)$ has at least two Leray-Hopf solutions. 
\end{Theorem}

When $0<\alpha=\beta<\frac12$, we automatically have the same result for system (\ref{sys-3}).

\begin{Remark}
In view of Theorem \ref{thm-unique-gmhd} and Theorem \ref{thm-nonunique2}, the value of $\frac12$ for $\alpha$ and $\beta$ is a sharp threshold to separate uniqueness from non-uniqueness result for system (\ref{sys-3}). On the other hand, for (\ref{sys-4}) the additional conditions of $\alpha\leq \beta$ and $3\beta-\alpha<1$ leave some gap where neither uniqueness nor non-uniqueness is known to hold.
\end{Remark}



\medskip

\subsection{Organisation of the paper}
We provide an outline of the rest of the paper. 
\begin{itemize}
\item
Section \ref{sec-pre} introduces notations and definitions of solutions for dyadic systems. 
\item
Section \ref{sec-LH} is devoted to a proof of Theorem \ref{thm-existence} on the existence of Leray-Hopf solutions. \item
Section \ref{sec-weak-strong} addresses the weak-strong uniqueness and uniqueness of Leray-Hopf solution for system (\ref{gen1})-(\ref{gen2}) with $\theta\leq 2$.
\item
In Section \ref{sec-non-unique} we construct non-unique Leray-Hopf solutions for (\ref{gen1})-(\ref{gen2}) with $\theta>2$.
\item
Section \ref{sec-gmhd-dyadic} outlines constructions to establish firstly conditions for uniqueness and secondly conditions for non-uniqueness of Leray-Hopf solutions to the dyadic model for MHD with fractional diffusions.
\end{itemize}

\bigskip

\section{Notations and notion of solutions}
\label{sec-pre}
\subsection{Notations}

The space $l^2$ is endowed with the standard scalar product and norm,
\[(u,v):=\sum_{n=1}^\infty u_nv_n, \ \ \ |u|:=\sqrt{(u,u)}.\]
It is regarded as the energy space in this paper.
We use $H^s$ to represent the space of sequences 
equipped with the scaler product 
\[(u,v)_s:=\sum_{n=1}^\infty \lambda_n^{2s}u_nv_n\]
and norm
\[\|u\|_s:=\sqrt{(u,u)_s}.\]
Notice that $H^0=l^2$. 

\subsection{Notion of solutions}
In the following, we introduce the concept of solutions for dyadic systems. 
\begin{Definition}\label{def1}
A pair of $l^2$-valued functions $(a(t), b(t))$ defined on $[0,\infty)$ is said to be a weak solution of (\ref{gen1})-(\ref{gen2}) if $a_j$ and $b_j$ satisfy (\ref{gen1})-(\ref{gen2}) and $a_j, b_j\in C^1([0,\infty))$ for all $j\geq0$.
\end{Definition}


\begin{Definition}\label{def-LH}
A Leray-Hopf solution $(a(t), b(t))$ of (\ref{gen1})-(\ref{gen2}) on $[0,T)$ is a weak solution satisfying\\
 \[a_j, b_j\in L^\infty([0,T); l^2)\cap L^2([0,T); H^1), \ \ \forall \ \ j\geq 0,\]
and
\begin{equation}\notag
\begin{split}
&\|a(t)\|_{l^2}^2+\|b(t)\|_{l^2}^2+2\nu\int_{0}^t\|a(\tau)\|_{H^1}^2 \, d\tau+2\mu\int_{0}^t\|b(\tau)\|_{H^1}^2 \, d\tau\\
\leq& \|a(0)\|_{l^2}^2+\|b(0)\|_{l^2}^2
\end{split}
\end{equation}
for all $0\leq t< T$. 
\end{Definition}

\begin{Definition}\label{def-LH}
A Leray-Hopf solution $(a(t), b(t))$ of (\ref{sys-3}) on $[0,T)$ is a weak solution satisfying\\
 \[a_j, b_j\in L^\infty([0,T); l^2)\cap L^2([0,T); H^{\alpha}), \ \ \forall \ \ j\geq 0,\]
and
\begin{equation}\notag
\begin{split}
&\|a(t)\|_{l^2}^2+\|b(t)\|_{l^2}^2+2\nu\int_{0}^t\|a(\tau)\|_{H^{\alpha}}^2 \, d\tau+2\mu\int_{0}^t\|b(\tau)\|_{H^{\alpha}}^2 \, d\tau\\
\leq& \|a(0)\|_{l^2}^2+\|b(0)\|_{l^2}^2
\end{split}
\end{equation}
for all $0\leq t< T$. 
\end{Definition}

Weak solution and Leray-Hopf solution of other dyadic systems in the paper can be defined analogously. 

To reduce the number of parameters, we take $\nu=\mu=1$ in the rest of the paper since they do not affect the estimates or constructions.

\bigskip

\section{Existence of Leray-Hopf solutions}
\label{sec-LH}

In this section, we apply the Galerkin approximating approach to show the existence of Leray-Hopf solutions to (\ref{gen1})-(\ref{gen2}). Since the value of $\kappa_1$ and $\kappa_2$ does not play a role in the proof, without loss of generality, we set $\kappa_1=-\kappa_2=1$.
Fix any integer $N\geq 1$. Denote the sequences 
\[a^N(t)=\{a_j^N(t)\}_{j\geq 0}, \ \ b^N(t)=\{b_j^N(t)\}_{j\geq 0}, \ \ \mbox{with} \ \  
a_j^N= b_j^N\equiv 0, \ \ \ \forall \ \ j\geq N+1.\]
That is,
\begin{equation}\notag
\begin{split}
a^N(t)=&\left(a_0^N(t), a_1^N(t), a_2^N(t), ..., a_N^N(t), 0, 0, 0, ... \right),\\
b^N(t)=&\left(b_0^N(t), b_1^N(t), b_2^N(t), ..., b_N^N(t), 0, 0, 0, ... \right).
\end{split}
\end{equation}
Consider the truncated system for $(a^N(t), b^N(t))$, 
\begin{equation}\label{sys-truncated}
\begin{split}
\frac{d}{dt}a_j^N=&-\lambda_j^2a_j^N-\lambda_j^{\theta}a_j^Na_{j+1}^N+\lambda_{j-1}^{\theta}(a_{j-1}^N)^2
-\lambda_{j}^{\theta}b_j^Nb_{j+1}^N\\
&+\lambda_{j-1}^{\theta}(b_{j-1}^N)^2+f_j,  \ \ \ \ 0\leq j\leq N\\
\frac{d}{dt}b_j^N=&-\lambda_j^2b_j^N+\lambda_j^{\theta}a_j^Nb_{j+1}^N-\lambda_{j}^{\theta}b_j^Na_{j+1}^N, \ \ \ \  0\leq j\leq N\\
a_j^N(0)=&\ a_j^0, \ \ b_j^N(0)= b_j^0, \ \ \ 0\leq j\leq N. 
\end{split}
\end{equation}
By convention, $a_{-1}^N=b_{-1}^N=0$.

In the following, we proceed with the standard Galerkin approximating framework: (i) for any $N\geq 1$, there is a solution $(a^N(t), b^N(t))$ to (\ref{sys-truncated}) with $a^N(t)$ and $b^N(t)$ in the space $L^\infty(0,T; l^2)\cap L^2(0,T; H^1)$ and satisfying the corresponding energy inequality; (ii) we pass the sequence $\left\{(a^N(t), b^N(t))\right\}_{N\geq 1}$ (or a subsequence of it) to a limit $(a(t), b(t))$; (iii) the limit $(a(t), b(t))$ is shown to be a Leray-Hopf solution of (\ref{gen1})-(\ref{gen2}).

The integral form of (\ref{sys-truncated}) is 
\begin{equation}\label{eq-integral}
\begin{split}
a_j^N(t)=&\ a_j^0+\int_0^t  \left(-\lambda_j^2a_j^N(\tau)-\lambda_j^{\theta}a_j^N(\tau)a_{j+1}^N(\tau)+\lambda_{j-1}^{\theta}(a_{j-1}^N(\tau))^2\right. \\
&\left. -\lambda_{j}^{\theta}b_j^N(\tau)b_{j+1}^N(\tau)+\lambda_{j-1}^{\theta}(b_{j-1}^N(\tau))^2+f_j(\tau)\right)\,d\tau,\\
b_j^N(t)=&\ b_j^0+\int_0^t\left(-\lambda_j^2b_j^N(\tau)+\lambda_j^{\theta}a_j^N(\tau)b_{j+1}^N(\tau)-\lambda_{j}^{\theta}b_j^N(\tau)a_{j+1}^N(\tau)\right)\, d\tau,
\end{split}
\end{equation}
for $0\leq j\leq N$. Denote 
\begin{equation}\notag
\begin{split}
F^N(a^N, b^N, t)=&\left(F_0^N(a^N, b^N, t), F_1^N(a^N, b^N, t), ..., F_N^N(a^N, b^N, t)\right), \\
G^N(a^N, b^N)=&\left(G_0^N(a^N, b^N), G_1^N(a^N, b^N), ..., G_N^N(a^N, b^N)\right),
\end{split}
\end{equation}
with 
\begin{equation}\notag
\begin{split}
F_j^N(a^N, b^N, t)=&-\lambda_j^2a_j^N(t)-\lambda_j^{\theta}a_j^N(t)a_{j+1}^N(t)+\lambda_{j-1}^{\theta}(a_{j-1}^N(t))^2 \\
& -\lambda_{j}^{\theta}b_j^N(t)b_{j+1}^N(t)+\lambda_{j-1}^{\theta}(b_{j-1}^N(t))^2+f_j(t),\\
G_j^N(a^N, b^N)=& -\lambda_j^2b_j^N(t)+\lambda_j^{\theta}a_j^N(t)b_{j+1}^N(t)-\lambda_{j}^{\theta}b_j^N(t)a_{j+1}^N(t),
\end{split}
\end{equation}
for $0\leq j\leq N$. Denote $a^{0,N}=(a_0^0, a_1^0, ..., a_N^0)$, $b^{0,N}=(b_0^0, b_1^0, ..., b_N^0)$ and $f^N=(f_0, f_1, ..., f_N)$. Thus, system (\ref{eq-integral}) can be written as
\begin{equation}\label{integral2}
\begin{split}
a^N(t)=&\ a^{0,N}+\int_0^t F^N(a^N(\tau), b^N(\tau), \tau)\, d\tau,\\
b^N(t)=&\ b^{0,N}+\int_0^t G^N(a^N(\tau), b^N(\tau))\, d\tau.
\end{split}
\end{equation}
Denote the map 
\begin{equation}\notag
M_N(a^N, b^N)(t)=
\begin{pmatrix}
a^{0,N}+\int_0^t F^N(a^N(\tau), b^N(\tau), \tau)\, d\tau \\ 
b^{0,N}+\int_0^t G^N(a^N(\tau), b^N(\tau))\, d\tau
\end{pmatrix}
.
\end{equation}
Notice that there exists a constant $C_N$ depending on $N$ such that 
\begin{equation}\label{FN}
\left| F^N(a^N, b^N, t)\right|\leq C_N\left(|a^N|+|a^N|^2+|b^N|^2\right)+|f^N|,
\end{equation}
\begin{equation}\label{GN}
\left| G^N(a^N, b^N, t)\right|\leq C_N\left(|b^N|+|a^N|^2+|b^N|^2\right),
\end{equation}
and moreover
\begin{equation}\label{FN2}
\begin{split}
&\left| F^N(a^N, b^N, t)-F^N(\tilde a^N, \tilde b^N, t) \right|\\
\leq &\ C_N\left(1+|a^N|+|\tilde a^N|+|b^N|+|\tilde b^N|\right)\left(|a^N-\tilde a^N|+|b^N-\tilde b^N|\right),
\end{split}
\end{equation}
\begin{equation}\label{GN2}
\begin{split}
&\left| G^N(a^N, b^N)-G^N(\tilde a^N, \tilde b^N) \right|\\
\leq &\ C_N\left(1+|a^N|+|b^N|+|\tilde a^N|+|\tilde b^N|\right)\left(|a^N-\tilde a^N|+|b^N-\tilde b^N|\right).
\end{split}
\end{equation}
Choose 
\begin{equation}\label{radius}
R_N=2|a^{0,N}|+2|b^{0,N}|+2\int_0^T\left|f^N(t)\right|\, dt,
\end{equation}
and 
\begin{equation}\label{time}
t_{N, 1}=\frac{1}{2C_N(2R_N+1)}.
\end{equation}
Consider the map $M_N(a^N, b^N)$ on the following closed subset of the space of continuous functions $C([0,t_{N,1}]; \mathbb R^N)$ 
\begin{equation}\notag
B_N=\left\{(u,v)\in C([0,t_{N,1}]; \mathbb R^N)\times C([0,t_{N,1}]; \mathbb R^N): \|u\|_C\leq R_N, \|v\|_C\leq R_N.\right\}
\end{equation}
We claim that $M_N$ is a contraction mapping on $B_N$. Indeed, for any $(a^N, b^N)\in B_N$, it follows from (\ref{FN}), (\ref{radius}) and (\ref{time}) that
\begin{equation}\notag
\begin{split}
&\left| a^{0,N}+\int_0^t F^N(a^N(\tau), b^N(\tau), \tau)\, d\tau \right|\\
\leq & \left|a^{0,N}\right|+\int_0^t \left|F^N(a^N(\tau), b^N(\tau), \tau)\right|\, d\tau\\
\leq & \left|a^{0,N}\right|+ tC_N\left(R_N+2R_N^2\right)+\int_0^T|f^N(t)|\, dt\\
\leq & \ \frac12 R_N+\frac{1}{2C_N(2R_N+1)}  C_N\left(R_N+2R_N^2\right)\\
=&\ R_N;
\end{split}
\end{equation} 
and similarly, by (\ref{GN}), (\ref{radius}) and (\ref{time})
\begin{equation}\notag
\begin{split}
&\left| b^{0,N}+\int_0^t G^N(a^N(\tau), b^N(\tau))\, d\tau \right|\\
\leq & \left|b^{0,N}\right|+\int_0^t \left|G^N(a^N(\tau), b^N(\tau))\right|\, d\tau\\
\leq & \left|b^{0,N}\right|+ tC_N\left(R_N+2R_N^2\right)\\
\leq & \ \frac12 R_N+\frac{1}{2C_N(2R_N+1)}  C_N\left(R_N+2R_N^2\right)\\
=&\ R_N.
\end{split}
\end{equation} 
Thus, $M_N$ maps $B_N$ to itself. On the other hand, the property of contraction follows from (\ref{FN2}), (\ref{GN2}) and the choice of time $t_{N,1}$ in (\ref{time}). Therefore, system (\ref{integral2}) has a solution $(a^N(t), b^N(t))$ on $[0, t_{N,1}]$, and so does system (\ref{sys-truncated}). Next we show that the solution satisfies the energy inequality. Multiplying the first equation of (\ref{sys-truncated}) by $a_j^N$ and the second one by $b_j^N$, taking the sum for $0\leq j\leq N$ and integrating over $[0,t]$, we obtain 
\begin{equation}\label{LF-energy0}
\begin{split}
&\sum_{j=0}^N\left(a_j^N(t)^2+b_j^N(t)^2\right)+2\sum_{j=0}^N\int_0^t \lambda_j^2\left(a_j^N(\tau)^2+b_j^N(\tau)^2\right)\, d\tau\\
=&\sum_{j=0}^N\left(\left(a_j^{0,N}\right)^2+\left(b_j^{0,N}\right)^2\right)+2\sum_{j=0}^N\int_0^t f_j(\tau)a_j^N(\tau)\, d\tau.
\end{split}
\end{equation}
Applying the Cauchy-Schwarz inequality, we have 
\begin{equation}\notag
2\sum_{j=0}^N\int_0^t f_j(\tau)a_j^N(\tau)\, d\tau\leq 
\sum_{j=0}^N\int_0^t \lambda_j^{-2}f_j^2(\tau)\, d\tau + \sum_{j=0}^N\int_0^t \lambda_j^2a_j^N(\tau)^2\, d\tau.
\end{equation}
Hence, it follows from (\ref{LF-energy0})
\begin{equation}\label{LF-energy1}
\begin{split}
&\sum_{j=0}^N\left(a_j^N(t)^2+b_j^N(t)^2\right)+\sum_{j=0}^N\int_0^t \lambda_j^2\left(a_j^N(\tau)^2+b_j^N(\tau)^2\right)\, d\tau\\
\leq & \sum_{j=0}^N\left(\left(a_j^{0,N}\right)^2+\left(b_j^{0,N}\right)^2\right)+\sum_{j=0}^N\int_0^t \lambda_j^{-2}f_j^2(\tau)\, d\tau.
\end{split}
\end{equation}
We can iterate the process above to construct the solution on time intervals $[t_{N,1}, t_{N,2}]$, $[t_{N,2}, t_{N,3}]$, ...,  $[t_{N,k}, t_{N,k+1}]$, ..., and finally reach the time $T$.  Indeed, we observe from the energy inequality (\ref{LF-energy1}) that 
\begin{equation}\notag
\left|a^N(t_k)\right|+\left|b^N(t_k)\right|\leq |a^{0,N}|+|b^{0,N}|+k\left(\sum_{j=0}^N\int_0^t \lambda_j^{-2}f_j^2(\tau)\, d\tau\right)^{\frac12}.
\end{equation}
Hence, according to (\ref{time}) and (\ref{radius}), at the $k+1$-th iteration, we have
\begin{equation}\notag
\begin{split}
&t_{N,k+1}-t_{N,k}\\
\geq& \left[C_N\left(8|a^{0,N}|+8|b^{0,N}|+4 \int_0^T|f^N(t)|\,dt+2+4k\left(\sum_{j=0}^N\int_0^t \lambda_j^{-2}f_j^2(\tau)\, d\tau\right)^{\frac12} \right)\right]^{-1}\\
\gtrsim& \frac{1}{k}.
\end{split}
\end{equation}
Therefore, the sum $\sum_{k}(t_{N,k+1}-t_{N,k})$ diverges and will reach $T$ after a certain number of iterations. In conclusion, we obtain a solution $(a^N(t), b^N(t))$ of (\ref{sys-truncated}) on the interval $[0,T]$, which satisfies the energy inequality (\ref{LF-energy1}) for all $t\in[0,T]$.

The next step is to extract a limit from the sequence $\{(a^N(t), b^N(t))\}_{N\geq 1}$. In view of the integral equations in (\ref{eq-integral}) and the energy inequality (\ref{LF-energy1}), we know $a_j^N, b_j^N\in C^1[0,T]$ for any $0\leq j\leq N$ and $a^N, b^N\in L^\infty(0,T; l^2)\cap L^2(0,T; H^1)$ for any $N\geq 1$. As a consequence, there exists a subsequence $\{(a^{N_k}(t), b^{N_k}(t))\}_{k\geq 1}$ which converges to $(a(t), b(t))$ in $C[0,T]$ such that (by employing a diagonal process)
\[a_j^{N_k} \to a_j, \ \ b_j^{N_k} \to b_j, \ \ \mbox{in} \ \ C(0,T) \ \ \mbox{as} \ \ k\to\infty, \ \ \forall \ \ j\geq 0.\]

The last step is to show that the limit $(a(t), b(t))$ is a Leray-Hopf solution of (\ref{gen1})-(\ref{gen2}). Replacing $N$ by $N_k$ in (\ref{eq-integral}) and taking the limit $k\to\infty$, we see that $(a(t), b(t))$ satisfies the integral system
\begin{equation}\notag
\begin{split}
a_j(t)=&\ a_j^0+\int_0^t  \left(-\lambda_j^2a_j(\tau)-\lambda_j^{\theta}a_j(\tau)a_{j+1}(\tau)+\lambda_{j-1}^{\theta}(a_{j-1}(\tau))^2\right. \\
&\left. -\lambda_{j}^{\theta}b_j(\tau)b_{j+1}(\tau)+\lambda_{j-1}^{\theta}(b_{j-1}(\tau))^2+f_j(\tau)\right)\,d\tau,\\
b_j(t)=&\ b_j^0+\int_0^t\left(-\lambda_j^2b_j(\tau)+\lambda_j^{\theta}a_j(\tau)b_{j+1}(\tau)-\lambda_{j}^{\theta}b_j(\tau)a_{j+1}(\tau)\right)\, d\tau,
\end{split}
\end{equation}
for all $j\geq 0$. Hence, $(a(t), b(t))$ satisfies system (\ref{gen1})-(\ref{gen2}). Moreover,  $a_j, b_j\in C^1[0,T]$ for all $j\geq 0$. In addition, taking the limit in the energy inequality (\ref{LF-energy1}) yields 
\[a, b\in L^\infty(0,T; l^2)\cap L^2(0,T; H^1). \] 
Notice that $a^{N_k}\in L^\infty(0,T; l^2)\cap L^2(0,T; H^1)$ for all $k$ and $N_k$. Thus, the sequence  $\{a_j^{N_k}\}_{k\geq 1}$ converges weakly in $L^2(0,T)$ for any fixed $j\geq0$, and the limit coincides with $a_j$. Consequently, we have
\begin{equation}\label{limit-f}
\sum_{j=0}^\infty \int_0^t f_j(\tau) a_j^{N_k}(\tau)\, d\tau \to \sum_{j=0}^\infty \int_0^t f_j(\tau) a_j(\tau)\, d\tau, \ \ \mbox{as} \ \ k\to \infty, \ \ \forall \ \ t\in[0,T].
\end{equation}
Passing the limit in (\ref{LF-energy0}) and applying (\ref{limit-f}), it leads to the energy inequality satisfied by the limit $(a(t), b(t))$
\begin{equation}\notag
\begin{split}
&\sum_{j=0}^\infty\left(a_j(t)^2+b_j(t)^2\right)+2\sum_{j=0}^\infty\int_0^t \lambda_j^2\left(a_j(\tau)^2+b_j(\tau)^2\right)\, d\tau\\
\leq &\sum_{j=0}^\infty\left(\left(a_j^{0}\right)^2+\left(b_j^{0}\right)^2\right)+2\sum_{j=0}^\infty\int_0^t f_j(\tau)a_j(\tau)\, d\tau.
\end{split}
\end{equation}
It completes the proof of Theorem \ref{thm-existence}.

\bigskip

\section{Weak-strong uniqueness }
\label{sec-weak-strong}

In order to show the weak-strong uniqueness, a standard argument involving Gr\"onwall's inequality will be applied to the difference $(a(t)-u(t), b(t)-v(t))$ of the two solutions $(a(t), b(t))$ and $(u(t), v(t))$. 
\medskip

{\textbf{Proof of Theorem \ref{thm-weak-strong}:}}
As in the previous section, we set $\kappa_1=-\kappa_2=1$. We start with the energy balance through the $j$-th shell
\begin{equation}\label{diff-energy1}
\begin{split}
&\frac{d}{dt}\left((a_j-u_j)^2+(b_j-v_j)^2\right)\\
=&\left((2a_ja_j'+2b_jb_j')+(2u_ju_j'+2v_jv_j')\right)\\
&-\left((2a_ju_j'+2a_j'u_j)+(2b_jv_j'+2b_j'v_j)\right)
\end{split}
\end{equation}
and continue to estimate the four groups on the right hand side.
In view of equations (\ref{gen1})-(\ref{gen2}) satisfied by $(a_j, b_j)$ and $(u_j, v_j)$ respectively, we have
\begin{equation}\label{diff-energy2}
a_ja_j'+b_jb_j'=-\lambda_j^2(a_j^2+b_j^2)-\lambda_j^\theta(a_j^2+b_j^2)a_{j+1}+\lambda_{j-1}^\theta(a_{j-1}^2+b_{j-1}^2)a_j+f_ja_j,
\end{equation}
\begin{equation}\label{diff-energy3}
u_ju_j'+v_jv_j'=-\lambda_j^2(u_j^2+v_j^2)-\lambda_j^\theta(u_j^2+v_j^2)u_{j+1}+\lambda_{j-1}^\theta(u_{j-1}^2+v_{j-1}^2)u_j+f_ju_j,
\end{equation}
\begin{equation}\label{diff-energy4}
\begin{split}
(a_ju_j)'=& -2\lambda_j^2a_ju_j-\lambda_j^\theta a_ju_ja_{j+1}-\lambda_j^\theta a_ju_ju_{j+1}\\
&-\lambda_j^\theta b_ju_jb_{j+1} -\lambda_j^\theta a_jv_jv_{j+1}
+\lambda_{j-1}^\theta a_{j-1}^2u_j+\lambda_{j-1}^\theta b_{j-1}^2u_j\\
&+\lambda_{j-1}^\theta u_{j-1}^2a_j+\lambda_{j-1}^\theta v_{j-1}^2a_j+f_j(a_j+u_j),
\end{split}
\end{equation}
\begin{equation}\label{diff-energy5}
\begin{split}
(b_jv_j)'=& -2\lambda_j^2b_jv_j+\lambda_j^\theta a_jv_jb_{j+1}+\lambda_j^\theta b_ju_jv_{j+1}-\lambda_j^\theta b_jv_ja_{j+1}-\lambda_j^\theta b_jv_ju_{j+1}.
\end{split}
\end{equation}
Combining (\ref{diff-energy1})-(\ref{diff-energy5}) and grouping the terms appropriately gives
\begin{equation}\label{diff-energy6}
\begin{split}
&\frac{d}{dt}\left((a_j-u_j)^2+(b_j-v_j)^2\right)+2\lambda_j^2(a_j-u_j)^2+2\lambda_j^2(b_j-v_j)^2\\
=&-\lambda_j^\theta(a_j^2+b_j^2)a_{j+1}+\lambda_{j-1}^\theta(a_{j-1}^2+b_{j-1}^2)a_j\\
&-\lambda_j^\theta(u_j^2+v_j^2)u_{j+1}+\lambda_{j-1}^\theta(u_{j-1}^2+v_{j-1}^2)u_j\\
&+\left(-\lambda_j^\theta a_ju_ja_{j+1}-\lambda_j^\theta a_ju_ju_{j+1}+\lambda_{j-1}^\theta a_{j-1}^2u_j+\lambda_{j-1}^\theta u_{j-1}^2a_j\right)\\
&+\left(-\lambda_j^\theta b_ju_jb_{j+1}-\lambda_j^\theta a_jv_jv_{j+1}+\lambda_j^\theta a_jv_jb_{j+1}+\lambda_j^\theta b_ju_jv_{j+1}\right)\\
&+\left(\lambda_{j-1}^\theta b_{j-1}^2u_j+\lambda_{j-1}^\theta v_{j-1}^2a_j-\lambda_j^\theta b_jv_ja_{j+1}-\lambda_j^\theta b_jv_ju_{j+1}\right).
\end{split}
\end{equation}
We further rearrange the terms in the last three parentheses of (\ref{diff-energy6}) to create terms in differences, for instance, $a_j-u_j$ and $b_j-v_j$. Shifting the sub-index $j$ to $j+1$ in the last two terms of 
\[\left(-\lambda_j^\theta a_ju_ja_{j+1}-\lambda_j^\theta a_ju_ju_{j+1}+\lambda_{j-1}^\theta a_{j-1}^2u_j+\lambda_{j-1}^\theta u_{j-1}^2a_j\right),\]
we have
\begin{equation}\label{diff-energy7}
\begin{split}
&-\lambda_j^\theta a_ju_ja_{j+1}-\lambda_j^\theta a_ju_ju_{j+1}+\lambda_{j-1}^\theta a_{j-1}^2u_j+\lambda_{j-1}^\theta u_{j-1}^2a_j\\
=&-\lambda_j^\theta a_ju_ja_{j+1}-\lambda_j^\theta a_ju_ju_{j+1}+\lambda_{j}^\theta a_{j}^2u_{j+1}+\lambda_{j}^\theta u_{j}^2a_{j+1}\\
=&\ \lambda_j^\theta\left(a_ju_{j+1}(a_j-u_j)-u_ja_{j+1}(a_j-u_j)\right)\\
=&\ \lambda_j^\theta (a_j-u_j) \left(a_ju_{j+1}-a_ja_{j+1}+a_ja_{j+1}-u_ja_{j+1}\right)\\
=& -\lambda_j^\theta a_j(a_j-u_j) \left(a_{j+1}-u_{j+1}\right)+\lambda_j^\theta\left(a_j-u_j\right)^2a_{j+1}.
\end{split}
\end{equation}
Similarly, with a shift of sub-index in the first two terms of 
\[\left(\lambda_{j-1}^\theta b_{j-1}^2u_j+\lambda_{j-1}^\theta v_{j-1}^2a_j-\lambda_j^\theta b_jv_ja_{j+1}-\lambda_j^\theta b_jv_ju_{j+1}\right),\]
we obtain
\begin{equation}\label{diff-energy8}
\begin{split}
&\lambda_{j}^\theta b_{j}^2u_{j+1}+\lambda_{j}^\theta v_{j}^2a_{j+1}-\lambda_j^\theta b_jv_ja_{j+1}-\lambda_j^\theta b_jv_ju_{j+1}\\
=&\ \lambda_{j}^\theta b_{j}u_{j+1}(b_j-v_j)-\lambda_{j}^\theta v_{j}a_{j+1}(b_j-v_j)\\
=&\ \lambda_{j}^\theta (b_j-v_j)(b_{j}u_{j+1}-b_ja_{j+1}+b_ja_{j+1}-v_ja_{j+1})\\
=& -\lambda_{j}^\theta b_j(b_j-v_j)(a_{j+1}-u_{j+1})+\lambda_{j}^\theta (b_j-v_j)^2 a_{j+1}.
\end{split}
\end{equation}
We rearrange the terms of 
\[\left(-\lambda_j^\theta b_ju_jb_{j+1}-\lambda_j^\theta a_jv_jv_{j+1}+\lambda_j^\theta a_jv_jb_{j+1}+\lambda_j^\theta b_ju_jv_{j+1}\right)\]
as 
\begin{equation}\label{diff-energy9}
\begin{split}
&-\lambda_j^\theta b_ju_jb_{j+1}-\lambda_j^\theta a_jv_jv_{j+1}+\lambda_j^\theta a_jv_jb_{j+1}+\lambda_j^\theta b_ju_jv_{j+1}\\
=&\ \lambda_j^\theta a_jv_j(b_{j+1}-v_{j+1})-\lambda_j^\theta b_ju_j (b_{j+1}-v_{j+1})\\
=&\ \lambda_j^\theta (a_jv_j-a_jb_j+a_jb_j-b_ju_j)(b_{j+1}-v_{j+1})\\
=&-\lambda_j^\theta a_j(b_j-v_j)(b_{j+1}-v_{j+1})+\lambda_j^\theta b_j(a_j-u_j)(b_{j+1}-v_{j+1}).
\end{split}
\end{equation}

Since $(a(t), b(t))$ and $(u(t), v(t))$ are Leray-Hopf solutions, we have that the following two series with telescope sums vanish,
\begin{equation}\label{diff-energy10}
\begin{split}
\sum_{j=0}^{\infty}\int_0^t \left(-\lambda_j^\theta(a_j^2+b_j^2)a_{j+1}+\lambda_{j-1}^\theta(a_{j-1}^2+b_{j-1}^2)a_j\right)  d\tau=&\ 0, \\
\sum_{j=0}^{\infty}\int_0^t \left(-\lambda_j^\theta(u_j^2+v_j^2)u_{j+1}+\lambda_{j-1}^\theta(u_{j-1}^2+v_{j-1}^2)u_j\right)  d\tau=&\ 0. \\
\end{split}
\end{equation}
Integrating (\ref{diff-energy6}) over $[0,t]$, taking the sum for $j\geq 0$, using the fact $a_{-1}=b_{-1}=u_{-1}=v_{-1}=0$, shifting the sub-index in the terms with sub-index $j-1$, and applying (\ref{diff-energy7})-(\ref{diff-energy10}), we deduce
\begin{equation}\label{diff-energy11}
\begin{split}
&\sum_{j=0}^{\infty}\left((a_j(t)-u_j(t))^2+(b_j(t)-v_j(t))^2\right)\\
&+2\sum_{j=0}^{\infty}\lambda_j^2\int_0^t(a_j(\tau)-u_j(\tau))^2+(b_j(\tau)-v_j(\tau))^2d\tau\\
=&-\sum_{j=0}^{\infty}\int_0^t \lambda_j^\theta a_j(\tau)(a_j(\tau)-u_j(\tau)) \left(a_{j+1}(\tau)-u_{j+1}(\tau)\right) d\tau\\
&-\sum_{j=0}^{\infty}\int_0^t \lambda_j^\theta b_j(\tau)(b_j(\tau)-v_j(\tau)) \left(a_{j+1}(\tau)-u_{j+1}(\tau)\right) d\tau\\
&+\sum_{j=0}^{\infty}\int_0^t \lambda_j^\theta (a_j(\tau)-u_j(\tau))^2 a_{j+1}(\tau) d\tau\\
& +\sum_{j=0}^{\infty}\int_0^t \lambda_j^\theta (b_j(\tau)-v_j(\tau))^2 a_{j+1}(\tau) d\tau\\
&-\sum_{j=0}^{\infty}\int_0^t \lambda_j^\theta a_j(\tau)(b_j(\tau)-v_j(\tau)) \left(b_{j+1}(\tau)-v_{j+1}(\tau)\right) d\tau\\
&+\sum_{j=0}^{\infty}\int_0^t \lambda_j^\theta b_j(\tau)(a_j(\tau)-u_j(\tau)) \left(b_{j+1}(\tau)-v_{j+1}(\tau)\right) d\tau.
\end{split}
\end{equation}
We claim that the series on the right hand side of (\ref{diff-energy11}) are well-defined.  Indeed, since $(a(t), b(t))$ and $(u(t), v(t))$ are Leray-Hopf solutions, it is clear that
\begin{equation}\notag
\sum_{j=0}^{\infty}\lambda_j^2\int_0^t a_j^2(\tau)+b_j^2(\tau)d\tau<\infty, \ \ \sum_{j=0}^{\infty}\lambda_j^2\int_0^t u_j^2(\tau)+v_j^2(\tau)d\tau<\infty.
\end{equation}
As a consequence, applying the assumption (\ref{assump-a}), we infer
\begin{equation}\notag
\begin{split}
&\sum_{j=J}^{\infty}\int_0^t \left|\lambda_j^\theta a_j(\tau)(a_j(\tau)-u_j(\tau)) \left(a_{j+1}(\tau)-u_{j+1}(\tau)\right) \right|d\tau\\
\leq &\ C \sum_{j=J}^{\infty}\int_0^t \left|\lambda_j^2(a_j(\tau)-u_j(\tau)) \left(a_{j+1}(\tau)-u_{j+1}(\tau)\right) \right|d\tau\\
\leq &\ 4C \sum_{j=J}^{\infty}\lambda_j^2\int_0^t a_j^2(\tau)+b_j^2(\tau)+u_j^2(\tau)+v_j^2(\tau)d\tau\\
< &\ \infty
\end{split}
\end{equation}
for a constant $C$.
Other series can be shown to converge analogously. Next, we estimate these series starting from the $J$-th shell. We only need to show details for one of them, for instance, thanks to the assumption (\ref{assump-a})
\begin{equation}\notag
\begin{split}
&\sum_{j=J}^{\infty}\int_0^t \left|\lambda_j^\theta b_j(\tau)(b_j(\tau)-v_j(\tau)) \left(a_{j+1}(\tau)-u_{j+1}(\tau)\right) \right|d\tau\\
\leq &\ C_0 \sum_{j=J}^{\infty}\int_0^t \left|\lambda_j^2(b_j(\tau)-v_j(\tau)) \left(a_{j+1}(\tau)-u_{j+1}(\tau)\right) \right|d\tau\\
\leq &\ \frac{C_0}2 \sum_{j=J}^{\infty}\int_0^t \lambda_j^2(b_j(\tau)-v_j(\tau))^2d\tau+\frac{C_0}{2\lambda^2} \sum_{j=J}^{\infty}\int_0^t \lambda_{j+1}^2(a_{j+1}(\tau)-u_{j+1}(\tau))^2d\tau.
\end{split}
\end{equation}
Similarly, the other series have the estimates
\begin{equation}\notag
\begin{split}
&\sum_{j=J}^{\infty}\int_0^t \left|\lambda_j^\theta a_j(\tau)(a_j(\tau)-u_j(\tau)) \left(a_{j+1}(\tau)-u_{j+1}(\tau)\right) \right|d\tau\\
\leq &\ \frac{C_0}2\sum_{j=J}^{\infty}\int_0^t \lambda_j^2(a_j(\tau)-u_j(\tau))^2d\tau+\frac{C_0}{2\lambda^2} \sum_{j=J}^{\infty}\int_0^t \lambda_{j+1}^2(a_{j+1}(\tau)-u_{j+1}(\tau))^2d\tau,\\
&\sum_{j=J}^{\infty}\int_0^t \left|\lambda_j^\theta a_j(\tau)(b_j(\tau)-v_j(\tau)) \left(b_{j+1}(\tau)-v_{j+1}(\tau)\right) \right|d\tau\\
\leq &\ \frac{C_0}2\sum_{j=J}^{\infty}\int_0^t \lambda_j^2(b_j(\tau)-v_j(\tau))^2d\tau+\frac{C_0}{2\lambda^2} \sum_{j=J}^{\infty}\int_0^t \lambda_{j+1}^2(b_{j+1}(\tau)-v_{j+1}(\tau))^2d\tau,\\
&\sum_{j=J}^{\infty}\int_0^t \left|\lambda_j^\theta b_j(\tau)(a_j(\tau)-u_j(\tau)) \left(b_{j+1}(\tau)-v_{j+1}(\tau)\right) \right|d\tau\\
\leq &\ \frac{C_0}2\sum_{j=J}^{\infty}\int_0^t \lambda_j^2(a_j(\tau)-u_j(\tau))^2d\tau+\frac{C_0}{2\lambda^2} \sum_{j=J}^{\infty}\int_0^t \lambda_{j+1}^2(b_{j+1}(\tau)-v_{j+1}(\tau))^2d\tau,\\
&\sum_{j=J}^{\infty}\int_0^t \left|\lambda_j^\theta (a_j(\tau)-u_j(\tau))^2 a_{j+1}(\tau)\right|d\tau
\leq  C_0\lambda^{2-\theta}\sum_{j=J}^{\infty}\int_0^t \lambda_j^2(a_j(\tau)-u_j(\tau))^2d\tau,\\
&\sum_{j=J}^{\infty}\int_0^t \left|\lambda_j^\theta (b_j(\tau)-v_j(\tau))^2 a_{j+1}(\tau)\right|d\tau
\leq  C_0\lambda^{2-\theta}\sum_{j=J}^{\infty}\int_0^t \lambda_j^2(b_j(\tau)-v_j(\tau))^2d\tau.
\end{split}
\end{equation}
Combining the estimates above and (\ref{diff-energy11}), we obtain
\begin{equation}\label{diff-energy12}
\begin{split}
&\sum_{j=0}^{\infty}\left((a_j(t)-u_j(t))^2+(b_j(t)-v_j(t))^2\right)\\
&+2\sum_{j=0}^{\infty}\lambda_j^2\int_0^t(a_j(\tau)-u_j(\tau))^2+(b_j(\tau)-v_j(\tau))^2d\tau\\
\leq & C_1\sum_{j=0}^J\int_0^t (a_j(\tau)-u_j(\tau))^2+(b_j(\tau)-v_j(\tau))^2d\tau\\
&+C_0\left(1+\lambda^{2-\theta}+\lambda^{-2}\right) \sum_{j=J}^{\infty}\lambda_j^2 \int_0^t (a_j(\tau)-u_j(\tau))^2+(b_j(\tau)-v_j(\tau))^2d\tau
\end{split}
\end{equation}
where the constant $C_1$ is given by
\begin{equation}\notag
C_1=32\lambda_J^\theta \sup_{0\leq j\leq J+1}\left(\|a_j\|_{C}+\|b_j\|_{C}\right)\leq 32\lambda_J^\theta\left(\|a^0\|_{l^2}+\|b^0\|_{l^2}\right).
\end{equation}
We take $C_0$ such that $C_0\left(1+\lambda^{2-\theta}+\lambda^{-2}\right) \leq 2$. Hence, it follows from (\ref{diff-energy12}) that
\begin{equation}\notag
\begin{split}
&\sum_{j=0}^{\infty}\left((a_j(t)-u_j(t))^2+(b_j(t)-v_j(t))^2\right)\\
\leq & C_1\sum_{j=0}^\infty\int_0^t (a_j(\tau)-u_j(\tau))^2+(b_j(\tau)-v_j(\tau))^2d\tau.
\end{split}
\end{equation}
Therefore, Gr\"onwall's inequality implies that 
\[a_j\equiv u_j, \ \ b_j\equiv v_j, \ \ \forall \ \ j\geq 0.\]

\cbdu

{\textbf{Proof of Theorem \ref{thm-unique1}:}} Since $0<\theta\leq 2$, for any Leray-Hopf solution $(a(t), b(t))$, there exists $J>0$ such that
\[|a_j(t)|\leq C_0 \lambda_j^{2-\theta}, \ \ |b_j(t)|\leq C_0 \lambda_j^{2-\theta}, \ \ \forall \ \ j\geq J. \]
 That is, assumption (\ref{assump-a}) is satisfied and hence uniqueness follows.
\cbdu

\bigskip

\section{Non-uniqueness of Leray-Hopf solutions for $\theta>2$}
\label{sec-non-unique}

We prove Theorem \ref{thm-nonunique1} in this section.
We adapt the construction scheme for the dyadic NSE in \cite{FK} in order to construct a solution $(a(t), b(t))$ of (\ref{gen1})-(\ref{gen2}) with zero initial data such that both $a(t)$ and $b(t)$ are non-vanishing. We first present the proof for the special case $\kappa_1=-\kappa_2=1$ and then point out modifications to prove other cases when changing the signs of $\kappa_1$ and $\kappa_2$.

Fix $T=\frac{1}{\lambda^2-1}$. Define 
\[t_j=\lambda_j^{-2}T, \ \ \ j\geq 0.\]
We note
\[t_{j-1}-t_j=\lambda_j^{-2}, \ \ \ j\geq 1,\]
\[(0,T)=\cup_{j=1}^\infty [t_j, t_{j-1}).\]
For $p,q\in C^\infty_0(0,1)$ and constant $\rho>\lambda^{\theta}$, we construct $a_j$ and $b_j$ as follows,
\begin{equation}\label{construct-a}
a_j(t)=
\begin{cases}
0, \ \ \ \ \ \ \ \ \ \ \ \ \ \ \ \ \ \ \ \ \ \ \ \ \ \ \ \ \ \ \ \ \ \ \ \ \ \ \ \ t<t_{j+1},\\
\lambda_{j+1}^{2-\theta} p\left(\lambda_{j+1}^2(t-t_{j+1})\right), \ \ \ \ \ \ \ t_{j+1}<t<t_j,\\
-\lambda_{j}^{2-\theta} q\left(\lambda_{j}^2(t-t_{j})\right), \ \ \ \ \ \ \ \ \ \ \ \ \ \ t_{j}<t<t_{j-1},\\
0, \ \ \ \ \ \ \ \ \ \ \ \ \ \ \ \ \ \ \ \ \ \ \ \ \ \ \ \ \ \ \ \ \ \ \ \ \ \ \ \ t>t_{j-1}.
\end{cases}
\end{equation}
\begin{equation}\label{construct-b}
b_j(t)=
\begin{cases}
0, \ \ \ \ \ \ \ \ \ \ \ \ \ \ \ \ \ \ \ \ \ \ \ \ \ \ \ \ \ \ \ \ \ \ \ \ \ \ \ \ \ \ \ t<t_{j+1},\\
\rho^{-j-1} h_1\left(\lambda_{j+1}^2(t-t_{j+1})\right), \ \ \ \ \ \ \ t_{j+1}<t<t_j,\\
\rho^{-j} h_2\left(\lambda_{j}^2(t-t_{j})\right), \ \ \ \ \ \ \ \ \ \ \ \ \ \ \ \ \ \ \ t_{j}<t<t_{j-1},\\
\rho^{-j+1} h_3\left(\lambda_{j-1}^2(t-t_{j-1})\right), \ \ \ \ \ \ \ t_{j-1}<t<t_{j-2},\\
\rho^{-j+1} h_3(1) e^{-\lambda_j^2(t-t_{j-2})}, \ \ \ \ \ \ \ \ \ \ \ \ \ \ \ t>t_{j-2},
\end{cases}
\end{equation}
such that $h_1, h_2$ and $h_3$ satisfy the ODE system on $[0,1]$
\begin{subequations}
\begin{align}
\frac{d}{dt} h_1+\left(\lambda^{-2}-\lambda^{-\theta}q\right) h_1-\lambda^{-\theta}ph_2=&\ 0,  \label{ode1} \\
\frac{d}{dt} h_2+h_2+q h_3=&\ 0, \label{ode2}\\
\frac{d}{dt} h_3+\lambda^2 h_3=&\ 0,  \label{ode3}\\
h_1(0)=0, \ \ h_2(0)=c_0, \ \ h_3(0)=&\ d_0.  \label{ode4}
\end{align}
\end{subequations}
In addition, we assume 
\begin{equation}\label{condition-h1}
h_1(1)= \rho c_0, \ \ \ h_2(1)= \rho d_0.
\end{equation}
With $(a_j, b_j)$ constructed in (\ref{construct-a})-(\ref{construct-b}), we define the forcing by 
\begin{equation}\label{construct-f}
f_j= \frac{d}{dt}a_j +\lambda_j^2 a_j
 + \lambda_j^{\theta}a_ja_{j+1}
+ \lambda_{j}^{\theta}b_jb_{j+1} -\lambda_{j-1}^{\theta}a_{j-1}^2-\lambda_{j-1}^{\theta}b_{j-1}^2
\end{equation}
for all $j\geq 0$. 

\begin{Lemma}\label{le-construct1}
Let $a_j$ and $b_j$ be constructed as in (\ref{construct-a})-(\ref{construct-b}). Then, the following properties hold:\\
(i) $a_j\in C^\infty_0(0,T)$ for all $j\geq 0$;\\
(ii) $b_j$ are piecewise smooth and $b_j \in H^1(0,T)$ for all $j\geq 0$;\\
(iii) \[a_j(0)=b_j(0)=0, \ \ \ \forall \ \ j\geq0;\] \\
(iv) \[ a_j(t)= O(\lambda_j^{2-\theta}), \ \ a_j'(t)= O(\lambda_j^{4-\theta}), \ \ b_j(t)= O(\rho^{-j}), \ \ \ j\to\infty.\]
\end{Lemma}
\pf
Since $p,q\in C^\infty_0(0,T)$, we only need to verify the values of the functions at $t_{j+1}$, $t_j$, $t_{j-1}$ and $t_{j-2}$.  The functions $b_j$ are piecewise smooth and continuous at these times; hence, $b_j \in H^1(0,T)$. It is obvious to see (iii) and (iv) from (\ref{construct-a})-(\ref{construct-b}).

\cbdu

\begin{Lemma}\label{le-construct2}
The functions $a_j$ and $b_j$ defined in (\ref{construct-a})-(\ref{construct-b}) satisfy 
\begin{equation}\notag
\frac{d}{dt}b_j +\lambda_j^2 b_j- \lambda_j^{\theta}a_jb_{j+1}+\lambda_{j}^{\theta}b_ja_{j+1}=0
\end{equation}
for all $j\geq 0$.
\end{Lemma}
\pf
Denote 
\[A_j(t)=\frac{d}{dt}b_j(t) +\lambda_j^2 b_j(t)- \lambda_j^{\theta}a_j(t)b_{j+1}(t)+\lambda_{j}^{\theta}b_j(t)a_{j+1}(t).\]
For $t<t_{j+1}$, we see $a_j(t)=b_j(t)=0$ from (\ref{construct-a}) and (\ref{construct-b}), and hence $A_j(t)=0$. 

For $t_{j+1}<t<t_j$, we denote $\tau=\lambda_{j+1}^2(t-t_{j+1})\in (0,1)$. It follows from (\ref{construct-a})-(\ref{construct-b}) that 
\begin{equation}\notag
\begin{split}
A_j(t)=&\ \rho^{-j-1}\lambda_{j+1}^2 \frac{d}{d\tau} h_1(\tau)+\rho^{-j-1}\lambda_j^2 h_1(\tau)\\
&-\rho^{-j-1}\lambda_{j+1}^2 \lambda^{-\theta}p(\tau) h_2(\tau)-\rho^{-j-1}\lambda_{j+1}^2 \lambda^{-\theta} q(\tau) h_1(\tau)\\
=&\ \rho^{-j-1}\lambda_{j+1}^2\left(\frac{d}{d\tau} h_1(\tau)+\lambda^{-2} h_1(\tau)- \lambda^{-\theta}q(\tau) h_1(\tau)- \lambda^{-\theta} p(\tau) h_2(\tau)\right)\\
=&\ 0
\end{split}
\end{equation}
thanks to (\ref{ode1}).

For $t_j<t<t_{j-1}$, we denote $\tau=\lambda_{j}^2(t-t_{j})\in (0,1)$. We note $a_{j+1}(t)=0$ by (\ref{construct-a}). Moreover, we have
\begin{equation}\notag
\begin{split}
A_j(t)=&\ \rho^{-j} \lambda_j^2 \frac{d}{d\tau} h_2(\tau)+\rho^{-j} \lambda_j^2 h_2(\tau)
+\rho^{-j} \lambda_j^{2} q(\tau) h_3(\tau)\\
=&\ \rho^{-j} \lambda_j^2\left(\frac{d}{d\tau} h_2(\tau)+ h_2(\tau)+ q(\tau)h_3(\tau)\right)\\
=&\ 0
\end{split}
\end{equation}
where we applied (\ref{ode2}). 

For $t_{j-1}<t<t_{j-2}$, we denote $\tau=\lambda_{j-1}^2(t-t_{j-1})\in (0,1)$. On this interval, we have $a_j(t)=a_{j+1}(t)=0$, and by (\ref{construct-b}) and (\ref{ode3})
\begin{equation}\notag
A_j(t)= \rho^{-j+1} \lambda_{j-1}^2\frac{d}{d\tau} h_3(\tau)+ \rho^{-j+1} \lambda_{j}^2 h_3(\tau)=0.
\end{equation}

For $t>t_{j-2}$, we note $a_j(t)=a_{j+1}(t)=0$, and 
\begin{equation}\notag
A_j(t)=-\rho^{-j+1}\lambda_j^2h_3(1) e^{-\lambda_j^2(t-t_{j-2})}+\rho^{-j+1}\lambda_j^2h_3(1) e^{-\lambda_j^2(t-t_{j-2})}=0.
\end{equation}

\cbdu

\begin{Lemma}\label{le-construct3}
The forcing $f=\{f_j(t)\}_{j\geq 0}$ constructed in (\ref{construct-f}) satisfies
\begin{equation}\notag
\sum_{j=0}^\infty \lambda_j^{-2}\int_0^T f_j^2(t)\, dt <\infty.
\end{equation}
\end{Lemma}
\pf
It follows from (\ref{construct-a})-(\ref{construct-b}), (\ref{construct-f}) and straightforward computations that
\begin{equation}\notag
f_j(t)=
\begin{cases}
0, \ \ \ \ \ \ \ \ \ t<t_{j+1}, \\
O(\lambda_j^{4-\theta}),  \ \ \ \ \ \ t_{j+1}<t<t_{j-2},\\
O(\lambda_j^{-\theta}),  \ \ \ \ \ \ t>t_{j-2}.
\end{cases}
\end{equation}
Hence,
\begin{equation}\notag
\begin{split}
\lambda_j^{-2}\int_0^T f_j^2(t)\, dt=&\ \lambda_j^{-2}\int_{t_{j+1}}^{t_{j-2}} O(\lambda_j^{8-2\theta})\, dt+
\lambda_j^{-2}\int_{t_{j-2}}^{T} O(\lambda_j^{-2\theta})\, dt\\
=&\ O(\lambda_j^{4-2\theta})+O(\lambda_j^{-2-2\theta}).
\end{split}
\end{equation}
Since $4-2\theta<0$ for $\theta>2$, it is clear that
\begin{equation}\notag
\sum_{j=0}^\infty \lambda_j^{-2}\int_0^T f_j^2(t)\, dt \leq \sum_{j=0}^\infty \left(O(\lambda_j^{4-2\theta})+O(\lambda_j^{-2-2\theta})\right)<\infty.
\end{equation}

\cbdu

\begin{Lemma}\label{le-construct4}
There exist functions $p,q\in C^\infty_0(0,1)$ and constants $c_0$ and $d_0$ with $c_0^2+d_0^2\neq 0$ such that there exists a unique solution $h=(h_1, h_2, h_3)$ of system (\ref{ode1})-(\ref{ode4}) satisfying (\ref{condition-h1}) and 
$h\in C^\infty([0,1]; \mathbb R^3)$.
\end{Lemma}
\pf
It is obvious from (\ref{ode3}) and the initial data $h_3(0)=d_0$ that
\[h_3(t)= d_0 e^{-\lambda^2 t}.\] 
It then follows from (\ref{ode2}) and $h_2(0)=c_0$ that 
\begin{equation}\notag
h_2(t)= c_0 e^{-t} - \int_0^t e^{s-t} q(s) h_3(s) \, ds = c_0 e^{-t} - d_0 e^{-t}\int_0^t e^{(1-\lambda^2)s} q(s) \, ds.
\end{equation}
Since $h_2(1)=\rho d_0$, we have the constraint
\begin{equation}\label{constraint1}
 c_0-d_0\int_0^1 e^{(1-\lambda^2)s} q(s) \, ds=e\rho d_0.
\end{equation}
In the end, we solve (\ref{ode1}) with $h_1(0)=0$ as
\begin{equation}\notag
h_1(t)=  \int_0^t  e^{-\int_s^t (\lambda^{-2}-\lambda^{-\theta}q(\tau))\, d\tau} \lambda^{-\theta} p(s) h_2(s)\, ds.
\end{equation}
The assumption $h_1(1)=\rho c_0$ gives another constraint, 
\begin{equation}\label{constraint2}
\int_0^1  e^{-\int_s^1 (\lambda^{-2}-\lambda^{-\theta}q(\tau))\, d\tau} \lambda^{-\theta} p(s) h_2(s)\, ds=\rho c_0.
\end{equation}
We note that in the case of constant $p$ and $q$, equations (\ref{constraint1})-(\ref{constraint2}) have a unique solution $(c_0,d_0)$. Thus, by a continuity argument, we know that there exist functions $p,q \in C^\infty_0(0,1)$ such that (\ref{constraint1}) and (\ref{constraint2}) are satisfied for some constants $c_0, d_0$ with $c_0^2+d_0^2\neq 0$. Since $p,q \in C^\infty_0(0,1)$, it is clear that $h_1,h_2, h_3\in C^\infty (0,1)$.

\cbdu

{\textbf{Proof of Theorem \ref{thm-nonunique1}:}}
Let $a=(a_j)_{j\geq 0}$ and $b=(b_j)_{j\geq 0}$ be constructed as in (\ref{construct-a})-(\ref{construct-b}). According to  Lemma \ref{le-construct2}, we have shown $(a,b)$ satisfies the model (\ref{gen1})-(\ref{gen2}) with $\kappa_1=-\kappa_2=1$ and with forcing $f_j$ defined in (\ref{construct-f}). It is shown in Lemma \ref{le-construct1} that $b\in H^1(0,T)$. We are left to show that $a\in l^2\cap H^1$ and $(a,b)$ satisfies the energy estimate. 

Since $\theta>2$, 
\begin{equation}\notag
a_j(t)=
\begin{cases}
0, \ \ \ \ \ \ \ \ \ \ \ \ \ \ \ \ \ \ \ \ \ t<t_{j+1},\\
O(\lambda_j^{2-\theta}),  \ \ \ \ \ \ \ t_{j+1}<t<t_{j-1},\\
0,  \ \ \ \ \ \ \ \ \ \ \ \ \  t>t_{j-1},\\
\end{cases}
\end{equation}
which implies 
\begin{equation}\notag
\sup_{t\in[0,T]}\sum_{j\geq 0}^\infty a_j^2(t)<\infty.
\end{equation}
Hence, we have $a\in l^2$. 

Notice 
\[t_{j-1}-t_j=\lambda_j^{-2}, \ \ t_j-t_{j+1}=\lambda_{j+1}^{-2}, \ \ T-t_{j-1} =\frac{1}{\lambda^2-1}\left(1-\frac{1}{\lambda_{j-1}^2}\right)<\frac{1}{\lambda^2-1}.\]
As a consequence, we have
\begin{equation}\notag
\begin{split}
\int_0^T a_j^2(t) \, dt=&\int_{t_{j+1}}^{t_{j-1}}a_j^2(t) \, dt
=\int_{t_{j+1}}^{t_{j-1}}O(\lambda_j^{4-2\theta}) \, dt
= O(\lambda_j^{2-2\theta}).
\end{split}
\end{equation}
Hence,
\begin{equation}\notag
\sum_{j=0}^\infty \int_0^T\lambda_j^2 a_j^2(t) \, dt
=\sum_{j=0}^\infty O(\lambda_j^{4-2\theta})<\infty
\end{equation}
provided $\theta>2$. That is, $a\in H^1(0,T)$. 

Next, we show that $(a,b)$ satisfies the energy identity. Since $(a(t), b(t))$ is a solution of (\ref{gen1})-(\ref{gen2}), it follows
\begin{equation}\notag
\begin{split}
\frac12\frac{d}{dt} \left(a_j^2(t)+b_j^2(t)\right)=&- \lambda_j^2 a_j^2- \lambda_j^2b_j^2-\lambda_j^\theta a_j^2a_{j+1}
+\lambda_{j-1}^\theta a_{j-1}^2 a_j\\
&-\lambda_j^\theta b_j^2 a_{j+1}+\lambda_{j-1}^\theta b_{j-1}^2 a_j+f_ja_j;
\end{split}
\end{equation}
and hence
\begin{equation}\label{energy-ajbj}
\begin{split}
&\frac{1}{2}\left(a_j^2(t)+b_j^2(t)\right)-\frac{1}{2}\left(a_j^2(0)+b_j^2(0)\right)\\
=&- \int_0^t\lambda_j^2 a_j^2\, dt- \int_0^t\lambda_j^2b_j^2\, dt-\int_0^t\lambda_j^\theta a_j^2a_{j+1}\, dt
+\int_0^t\lambda_{j-1}^\theta a_{j-1}^2 a_j\, dt\\
&-\int_0^t\lambda_j^\theta b_j^2 a_{j+1}\, dt+\int_0^t\lambda_{j-1}^\theta b_{j-1}^2 a_j\, dt+\int_0^t f_ja_j\, dt.
\end{split}
\end{equation}
Again, notice 
\begin{equation}\notag
a_j(t)=
\begin{cases}
0, \ \ \ \ \ \ \ \ \ \ \ \ \ \ \ \ \ \ \ \ \ t<t_{j+1},\\
O(\lambda_j^{2-\theta}),  \ \ \ \ \ \ \ t_{j+1}<t<t_{j-1},\\
0,  \ \ \ \ \ \ \ \ \ \ \ \ \  t>t_{j-1},\\
\end{cases}
\end{equation}
and 
\begin{equation}\notag
a_{j+1}(t)=
\begin{cases}
0, \ \ \ \ \ \ \ \ \ \ \ \ \ \ \ \ \ \ \ \ \ t<t_{j+2},\\
O(\lambda_{j+1}^{2-\theta}),  \ \ \ \ \ \ \ t_{j+2}<t<t_{j},\\
0,  \ \ \ \ \ \ \ \ \ \ \ \  t>t_{j}.\\
\end{cases}
\end{equation}
Thus, we have
\begin{equation}\notag
\begin{split}
\int_0^t\lambda_j^\theta \left|a_j^2a_{j+1}\right|\, dt=& \int_{t_{j+2}}^{t_{j-1}}\lambda_j^\theta O(\lambda_j^{6-3\theta})\, dt 
=  O(\lambda_j^{\theta-2+6-3\theta}) 
< \infty
\end{split}
\end{equation}
since $\theta>2$. Obviously, we also have
\begin{equation}\notag
\int_0^t\lambda_{j-1}^\theta \left|a_{j-1}^2a_{j}\right|\, dt<\infty.
\end{equation}
On the other hand, we note
\begin{equation}\notag
b_j(t)=
\begin{cases}
0, \ \ \ \ \ \ \ \ \ \ \ \ \ \ \ \ \ \ \ \ \ \ t<t_{j+1},\\
O(\rho^{-j-1}),  \ \ \ \ \ \ \ t_{j+1}<t<t_{j},\\
O(\rho^{-j}),  \ \ \ \ \ \ \ \ \ \ \ \ \  t_j<t<t_{j-1},\\
O(\rho^{-j+1}),  \ \ \ \ \ \ \ \ \ \ \ \  t>t_{j-1}.\\
\end{cases}
\end{equation}
Thus, we have
\begin{equation}\notag
\begin{split}
\int_0^t\lambda_j^\theta \left|b_j^2a_{j+1}\right|\, dt=& \int_{t_{j+1}}^{t_{j}}\lambda_j^\theta O(\rho^{-2j-2}\lambda_{j+1}^{2-\theta})\, dt+\int_{t_{j}}^{t_{j-1}}\lambda_j^\theta O(\rho^{-2j}\lambda_{j+1}^{-\theta})\, dt\\
&+ \int_{t_{j-1}}^{T}\lambda_j^\theta O(\rho^{-2j+2}\lambda_{j+1}^{-\theta})\, dt\\
=&\  O(\rho^{-2j-2})+O(\rho^{-2j}\lambda_{j}^{-2})+O(\rho^{-2j+2})\\
<&\ \infty.
\end{split}
\end{equation}
Similarly,
\begin{equation}\notag
\int_0^t\lambda_{j-1}^\theta \left|b_{j-1}^2a_{j}\right|\, dt<\infty.
\end{equation}
Therefore, we can take the sum of (\ref{energy-ajbj}) over $j\geq 0$ and obtain
\begin{equation}\notag
\begin{split}
&\frac{1}{2}\sum_{j=0}^\infty \left(a_j^2(t)+b_j^2(t)\right)-\frac{1}{2}\sum_{j=0}^\infty\left(a_j^2(0)+b_j^2(0)\right)\\
=&-\sum_{j=0}^\infty \int_0^t\lambda_j^2 a_j^2\, dt- \sum_{j=0}^\infty\int_0^t\lambda_j^2b_j^2\, dt
+\sum_{j=0}^\infty\int_0^t f_ja_j\, dt.
\end{split}
\end{equation}
Thus, we conclude $(a,b)$ is a Leray-Hopf solution of (\ref{gen1})-(\ref{gen2}) with zero initial data; however, $a\neq 0$ and $b\neq 0$. Non-uniqueness then follows. Indeed, for such forcing $f(t)$ as in (\ref{construct-f}), considering $b(t)\equiv 0$ in (\ref{gen1})-(\ref{gen2}), the forced dyadic model of the NSE has a solution $\tilde a(t)$. Hence, $(\tilde a(t), 0)$ is a trivial solution for the dyadic MHD model (\ref{gen1})-(\ref{gen2}).

\cbdu

\begin{Remark}
When $\kappa_1$ and $\kappa_2$ take different signs,
we can choose the same constructions for $(a_j(t), b_j(t))$ as in (\ref{construct-a}) and (\ref{construct-b}). The difference comes in the ODE system (\ref{ode1})-(\ref{ode4}) for the profile functions $h_1$, $h_2$ and $h_3$. For instance, if $\kappa_1=\kappa_2=1$ 
the functions $h_1$, $h_2$ and $h_3$ satisfy the following system 
\begin{subequations}
\begin{align}
\frac{d}{dt} h_1+\left(\lambda^{-2}+\lambda^{-\theta}q\right) h_1+\lambda^{-\theta}ph_2=&\ 0,  \label{ode1-fb} \\
\frac{d}{dt} h_2+h_2-q h_3=&\ 0, \label{ode2-fb}\\
\frac{d}{dt} h_3+\lambda^2 h_3=&\ 0,  \label{ode3-fb}\\
h_1(0)=0, \ \ h_2(0)=c_0, \ \ h_3(0)=&\ d_0,  \label{ode4-fb}
\end{align}
\end{subequations}
accompanied with 
\begin{equation}\label{condition-h1-fb}
h_1(1)= \rho c_0, \ \ \ h_2(1)= \rho d_0.
\end{equation}
We note that the structure of system (\ref{ode1-fb})-(\ref{ode4-fb}) remains similar to that of system (\ref{ode1})-(\ref{ode4}).
Thus in analogy with Lemma \ref{le-construct4}, it is not hard to show the existence of a solution $(h_1, h_2, h_3)$ to system (\ref{ode1-fb})-(\ref{ode4-fb}) satisfying (\ref{condition-h1-fb}). The rest analysis of Section \ref{sec-non-unique}  also holds for system (\ref{gen1})-(\ref{gen2}) with $\kappa_1=\kappa_2=1$.
\end{Remark}


\section{Uniqueness and non-uniqueness results for the dyadic MHD model with fractional Laplacians}
\label{sec-gmhd-dyadic}

\subsection{Uniqueness}
The weak-strong uniqueness stated in Theorem \ref{thm-weak-strong-gmhd} under assumption (\ref{assump-gmhd}) can be proved by following the steps described in Section \ref{sec-weak-strong}. We briefly present the main steps and emphasize why assumption (\ref{assump-gmhd}) is required for the uniqueness. 

Let $(a(t), b(t))$ and $(u(t), v(t))$ be two Leray-Hopf solutions of (\ref{sys-4}) with $(a(t), b(t))$ satisfying (\ref{assump-gmhd}). The difference of the two solutions satisfies the energy estimate
\begin{equation}\label{diff1-energy-gmhd}
\begin{split}
&\sum_{j=0}^{\infty}\left((a_j(t)-u_j(t))^2+(b_j(t)-v_j(t))^2\right)\\
&+2\sum_{j=0}^{\infty}\int_0^t\lambda_j^{2\alpha}(a_j(\tau)-u_j(\tau))^2+\lambda_j^{2\beta}(b_j(\tau)-v_j(\tau))^2d\tau\\
= &-\sum_{j=0}^{\infty}\int_0^t \lambda_j a_j(\tau)(a_j(\tau)-u_j(\tau)) \left(a_{j+1}(\tau)-u_{j+1}(\tau)\right) d\tau\\
&-\sum_{j=0}^{\infty}\int_0^t \lambda_j b_j(\tau)(b_j(\tau)-v_j(\tau)) \left(a_{j+1}(\tau)-u_{j+1}(\tau)\right) d\tau\\
&+\sum_{j=0}^{\infty}\int_0^t \lambda_j (a_j(\tau)-u_j(\tau))^2 a_{j+1}(\tau) d\tau
 +\sum_{j=0}^{\infty}\int_0^t \lambda_j (b_j(\tau)-v_j(\tau))^2 a_{j+1}(\tau) d\tau\\
&-\sum_{j=0}^{\infty}\int_0^t \lambda_j a_j(\tau)(b_j(\tau)-v_j(\tau)) \left(b_{j+1}(\tau)-v_{j+1}(\tau)\right) d\tau\\
&+\sum_{j=0}^{\infty}\int_0^t \lambda_j b_j(\tau)(a_j(\tau)-u_j(\tau)) \left(b_{j+1}(\tau)-v_{j+1}(\tau)\right) d\tau\\
\equiv&\ I_1+I_2+I_3+I_4+I_5+I_6.
\end{split}
\end{equation}
With $(a(t), b(t))$ and $(u(t), v(t))$ being Leray-Hopf solutions of (\ref{sys-4}), it holds
\begin{equation}\notag
\sum_{j=0}^{\infty}\int_0^t \left(\lambda_j^{2\alpha}a_j^2(\tau)+\lambda_j^{2\beta}b_j^2(\tau)\right)d\tau<\infty, \ \ \sum_{j=0}^{\infty}\int_0^t \left(\lambda_j^{2\alpha}u_j^2(\tau)+\lambda_j^{2\beta}v_j^2(\tau)\right)d\tau<\infty.
\end{equation}
Thus combining assumption (\ref{assump-gmhd}) we know that all of the series on the right hand side of (\ref{diff1-energy-gmhd}) are well-defined. 
Moreover, these series can be estimated in the following way by using (\ref{assump-gmhd}). For example, we estimate $I_2$ thanks to the condition $|b_j|\leq C_0\lambda_j^{\alpha+\beta-1}$ of (\ref{assump-gmhd}),
\begin{equation}\notag
\begin{split}
&\sum_{j=J}^{\infty}\int_0^t \lambda_j\left| b_j(\tau)(b_j(\tau)-v_j(\tau)) \left(a_{j+1}(\tau)-u_{j+1}(\tau)\right) \right|d\tau\\
\leq &\ C_0 \sum_{j=J}^{\infty}\int_0^t \left|\lambda_j^{\alpha+\beta}(b_j(\tau)-v_j(\tau)) \left(a_{j+1}(\tau)-u_{j+1}(\tau)\right) \right|d\tau\\
\leq &\ \frac{C_0}2 \sum_{j=J}^{\infty}\int_0^t \lambda_j^{2\beta}(b_j(\tau)-v_j(\tau))^2d\tau+\frac{C_0}{2\lambda^{2\alpha}} \sum_{j=J}^{\infty}\int_0^t \lambda_{j+1}^{2\alpha}(a_{j+1}(\tau)-u_{j+1}(\tau))^2d\tau.
\end{split}
\end{equation}
The term $I_6$ can be handled similarly. The condition $|a_j|\leq C_0\lambda_j^{2\alpha-1}$ is posed to estimate $I_1$ and $I_3$, and $|a_j|\leq C_0\lambda_j^{2\beta-1}$ is for $I_4$ and $I_5$. 
With the estimates, it follows from (\ref{diff1-energy-gmhd}) that
\begin{equation}\label{diff2-energy-gmhd}
\begin{split}
&\sum_{j=0}^{\infty}\left((a_j(t)-u_j(t))^2+(b_j(t)-v_j(t))^2\right)\\
&+2\sum_{j=0}^{\infty}\lambda_j^2\int_0^t(a_j(\tau)-u_j(\tau))^2+(b_j(\tau)-v_j(\tau))^2d\tau\\
\leq & C_1\sum_{j=0}^J\int_0^t (a_j(\tau)-u_j(\tau))^2+(b_j(\tau)-v_j(\tau))^2d\tau\\
&+C_0\left(1+\lambda^{2\alpha-1}+\lambda^{-2\alpha}\right) \sum_{j=J}^{\infty} \int_0^t \lambda_j^{2\alpha}(a_j(\tau)-u_j(\tau))^2d\tau\\
&+C_0\left(1+\lambda^{2\beta-1}+\lambda^{-2\beta}\right) \sum_{j=J}^{\infty} \int_0^t \lambda_j^{2\alpha}(b_j(\tau)-v_j(\tau))^2d\tau
\end{split}
\end{equation}
with a constant $C_1>0$. The constant $C_0$ can be chosen small enough such that
$C_0\left(1+\lambda^{2\alpha-1}+\lambda^{-2\alpha}\right) \leq 2$ and $C_0\left(1+\lambda^{2\beta-1}+\lambda^{-2\beta}\right) \leq 2$. Consequently, we have from (\ref{diff2-energy-gmhd}) that
\begin{equation}\notag
\begin{split}
&\sum_{j=0}^{\infty}\left((a_j(t)-u_j(t))^2+(b_j(t)-v_j(t))^2\right)\\
\leq&\  C_1\sum_{j=0}^\infty\int_0^t (a_j(\tau)-u_j(\tau))^2+(b_j(\tau)-v_j(\tau))^2d\tau.
\end{split}
\end{equation}
Gr\"onwall's inequality immediately implies that 
$a_j\equiv u_j$ and $b_j\equiv v_j$ for all $j\geq 0$.

\medskip

\subsection{Non-uniqueness}

The construction scheme to prove Theorem \ref{thm-nonunique2} is similar to that presented in Section \ref{sec-non-unique}. The main effort is to determine the scaling in constructing $a_j$ and $b_j$. To be complete, we specify the constructions as follows. 
For $T=\frac{1}{\lambda^{2\beta}-1}$, we take the partition 
\[t_j=\lambda_j^{-2\beta}T, \ \ \ j\geq 0,\]
with 
\[t_{j-1}-t_j=\lambda_j^{-2\beta}, \ \ \ j\geq 1; \ \ \ (0,T)=\cup_{j=1}^\infty [t_j, t_{j-1}).\]
For $p,q\in C^\infty_0(0,1)$ and constant $\rho>\lambda$, we choose $a_j$ and $b_j$ as
\begin{equation}\label{construct-a-gmhd}
a_j(t)=
\begin{cases}
0, \ \ \ \ \ \ \ \ \ \ \ \ \ \ \ \ \ \ \ \ \ \ \ \ \ \ \ \ \ \ \ \ \ \ \ \ \ \ \ \ t<t_{j+1},\\
\lambda_{j+1}^{2\beta-1} p\left(\lambda_{j+1}^{2\beta}(t-t_{j+1})\right), \ \ \ \ \ \ \ t_{j+1}<t<t_j,\\
-\lambda_{j}^{2\beta-1} q\left(\lambda_{j}^{2\beta}(t-t_{j})\right), \ \ \ \ \ \ \ \ \ \ \ \ \ \ t_{j}<t<t_{j-1},\\
0, \ \ \ \ \ \ \ \ \ \ \ \ \ \ \ \ \ \ \ \ \ \ \ \ \ \ \ \ \ \ \ \ \ \ \ \ \ \ \ \ t>t_{j-1}.
\end{cases}
\end{equation}
\begin{equation}\label{construct-b-gmhd}
b_j(t)=
\begin{cases}
0, \ \ \ \ \ \ \ \ \ \ \ \ \ \ \ \ \ \ \ \ \ \ \ \ \ \ \ \ \ \ \ \ \ \ \ \ \ \ \ \ \ \ \ t<t_{j+1},\\
\rho^{-j-1} h_1\left(\lambda_{j+1}^{2\beta}(t-t_{j+1})\right), \ \ \ \ \ \ \ t_{j+1}<t<t_j,\\
\rho^{-j} h_2\left(\lambda_{j}^{2\beta}(t-t_{j})\right), \ \ \ \ \ \ \ \ \ \ \ \ \ \ \ \ \ \ \ t_{j}<t<t_{j-1},\\
\rho^{-j+1} h_3\left(\lambda_{j-1}^{2\beta}(t-t_{j-1})\right), \ \ \ \ \ \ \ t_{j-1}<t<t_{j-2},\\
\rho^{-j+1} h_3(1) e^{-\lambda_j^{2\beta}(t-t_{j-2})}, \ \ \ \ \ \ \ \ \ \ \ \ \ \ \ t>t_{j-2},
\end{cases}
\end{equation}
where $h_1, h_2$ and $h_3$ are functions satisfying the following ODE system on $[0,1]$
\begin{subequations}
\begin{align}
\frac{d}{dt} h_1+\left(\lambda^{-2\beta}-\lambda^{-1}q\right) h_1-\lambda^{-1}ph_2=&\ 0,  \label{ode1-gmhd} \\
\frac{d}{dt} h_2+h_2+q h_3=&\ 0, \label{ode2-gmhd}\\
\frac{d}{dt} h_3+\lambda^{2\beta} h_3=&\ 0,  \label{ode3-gmhd}\\
h_1(0)=0, \ \ h_2(0)=c_0, \ \ h_3(0)=&\ d_0,  \label{ode4-gmhd}\\
h_1(1)= \rho c_0, \ \ \ h_2(1)=&\ \rho d_0. \label{ode5-gmhd}
\end{align}
\end{subequations}
We take the forcing $f_j$ as 
\begin{equation}\label{construct-f-gmhd}
f_j= \frac{d}{dt}a_j +\lambda_j^{2\alpha} a_j
 + \lambda_ja_ja_{j+1}
+ \lambda_{j}b_jb_{j+1} -\lambda_{j-1}a_{j-1}^2-\lambda_{j-1}b_{j-1}^2
\end{equation}
for all $j\geq 0$. 

For $f_j$ defined by (\ref{construct-f-gmhd}), we can show that $(a(t), b(t))$ with components constructed as in (\ref{construct-a-gmhd})-(\ref{construct-b-gmhd}) is a Leray-Hopf solution of system (\ref{sys-4}) with non-vanishing $b(t)$. We state the main ingredients to prove Theorem \ref{thm-nonunique2} in the following Lemmas, the proof of which are omitted. 

\begin{Lemma}\label{le-construct1-gmhd}
Let $\rho>\lambda>1$ and $0<\alpha\leq\beta<\frac12$.
The following properties hold for $a_j$ and $b_j$ as constructed in (\ref{construct-a-gmhd})-(\ref{construct-b-gmhd}):\\
(i) $a_j\in C^\infty_0(0,T)$ and $a_j \in H^\alpha(0,T)$ for all $j\geq 0$;\\
(ii) $b_j$ are piecewise smooth and $b_j \in H^\beta(0,T)$ for all $j\geq 0$;\\
(iii) \[a_j(0)=b_j(0)=0, \ \ \ \forall \ \ j\geq0;\] \\
(iv) \[ a_j(t)= O(\lambda_j^{2\beta-1}), \ \ a_j'(t)= O(\lambda_j^{4\beta-1}), \ \ b_j(t)= O(\rho^{-j}), \ \ \ j\to\infty.\]
\end{Lemma}

\begin{Lemma}\label{le-construct2-gmhd}
The functions $a_j$ and $b_j$ defined in (\ref{construct-a-gmhd})-(\ref{construct-b-gmhd}) satisfy system (\ref{sys-4}) with forcing $f_j$ defined by (\ref{construct-f-gmhd}).
\end{Lemma}

\begin{Lemma}\label{le-construct3-gmhd}
Let $\rho>\lambda>1$. Assume $0<\alpha\leq\beta<\frac12$ and $3\beta-\alpha<1$. The forcing $f_j$ defined by (\ref{construct-f-gmhd}) satisfies
\begin{equation}\notag
\sum_{j=0}^\infty \lambda_j^{-2\alpha}\int_0^T f_j^2(t)\, dt <\infty.
\end{equation}
\end{Lemma}

\begin{Lemma}\label{le-construct4-gmhd}
Let $\rho>\lambda>1$ and $0<\alpha\leq\beta<\frac12$. There exist functions $p,q\in C^\infty_0(0,1)$ and constants $c_0$, $d_0$ satisfying $c_0^2+d_0^2\neq 0$ such that there exists a unique solution $h=(h_1,h_2,h_3)\in C^\infty([0,1]; \mathbb R^3)$ of system (\ref{ode1-gmhd})-(\ref{ode5-gmhd}).
\end{Lemma}

\begin{Lemma}\label{le-construct5-gmhd}
Let $\rho>\lambda>1$ and $0<\alpha\leq\beta<\frac12$. Then $(a(t), b(t))$ satisfies the energy identity
\begin{equation}\notag
\begin{split}
&\frac{1}{2}\sum_{j=0}^\infty \left(a_j^2(t)+b_j^2(t)\right)-\frac{1}{2}\sum_{j=0}^\infty\left(a_j^2(0)+b_j^2(0)\right)\\
=&-\sum_{j=0}^\infty \int_0^t\lambda_j^{2\alpha} a_j^2\, dt- \sum_{j=0}^\infty\int_0^t\lambda_j^{2\beta}b_j^2\, dt
+\sum_{j=0}^\infty\int_0^t f_ja_j\, dt
\end{split}
\end{equation}
with $f_j$ defined by (\ref{construct-f-gmhd}).
\end{Lemma}

\bigskip

\section*{Acknowledgement}
M. Dai is partially supported by the NSF grants DMS-1815069 and DMS-2009422. S. Friedlander is partially supported by the NSF grant DMS-1613135. S. Friedlander is grateful to IAS for its hospitality in 2020-2021. M. Dai is also grateful to IAS for its hospitality in 2021-2022.

\bigskip


\end{document}